\documentclass [a4,11pt]{amsart}

\usepackage{amsmath,amssymb,amsfonts,enumerate,amsthm, amscd,}

\newcommand{\R}{\mathbb{R}}

\newcommand{\N}{\mathbb{N}}

\newcommand{\fm}{\mathfrak{m}}
\newcommand{\fp}{\mathfrak{p}}
\newcommand{\fq}{\mathfrak{q}}

\newtheorem{theorem}{Theorem}[section]
\newtheorem{lemma}[theorem]{Lemma}
\newtheorem{proposition}[theorem]{Proposition}
\newtheorem{corollary}[theorem]{Corollary}

\theoremstyle{definition}
\newtheorem{definition}[theorem]{Definition}

\theoremstyle{remark}
\newtheorem{remark}[theorem]{Remark}
\newtheorem{example}[theorem]{Example}

\theoremstyle{Definition and Notation}

\begin{document}
\bibliographystyle{amsplain}


\title[On Weakly von Neumann regular rings]{On Weakly von Neumann regular rings}

\author{Mohammed Kabbour}
\address{Mohammed Kabbour\\Department of Mathematics, Faculty of Science and Technology of Fez, Box 2202, University S.M. Ben Abdellah Fez, Morocco.
$$ E-mail\ address:\ mkabbour@gmail.com$$}

\author{Najib Mahdou}
\address{Najib Mahdou\\Department of Mathematics, Faculty of Science and Technology of Fez, Box 2202, University S.M. Ben Abdellah Fez, Morocco.
 $$E-mail\  address:\ mahdou@hotmail.com$$}



\begin{abstract}

In this paper, we define and study a particular case of von
Neumann regular notion called a weak von Neumann regular ring. It
shown that the polynomial ring $R[x]$ is weak von Neumann regular
if and only if $R$ has exactly two idempotent elements. We provide
necessary and sufficient conditions for $ R=A\propto E $ to be a
weak von Neumann ring. It is also shown that $I$ is a primary
ideal imply $R/I$ is a weak von Neumann regular ring.

\end{abstract}

\maketitle

 \begin{section} {\ Introduction}
 All rings considered in this paper are assumed to be commutative,
 and have identity element $\not=0;$ all modules are unital. A
 ring $R$ is reduced if its nilradical is zero. The following
 statements on a ring $R$ are equivalent:\begin{enumerate}
    \item Every finitely generated ideal of $R$ is principal and is
    generated by an idempotent.
    \item For each $x$ in $R,$ there is some $y$ in $R$ such that
    $x=x^2y.$
    \item $R$ is reduced 0-dimensional ring.
 \end{enumerate}
 A ring satisfying the equivalent conditions as above is said to
 be von Neumann regular. See for instance \cite{Ar,Gl,H,R}.
  \par In this article we study a new concept, close to the notion of
  von Neumann regular ring. More exactly we modify
  ( every finitely generated ideal $I\subseteq R$) by (every finitely
  generated ideal $I\subseteq J=Re\varsubsetneq R,$ where $e$ is an idempotent element of $R$) in the assertion (1).
  In the second section we give some results that allow us to study this
   notion, and it is containing some applications of such a notion.

\par A ring $R$ is called a coherent ring if every finitely generated ideal of
$R$ is finitely presented. We say that $R$ is coherent ring if and
only if $(0:a)$ is finitely generated ideal for every element $a$
of $R$ and the intersection of two finitely generated ideals of
$R$ is a finitely generated ideal of $R.$ Hence every von Neuman
regular ring is a coherent ring [ \cite{Gl}, p.47 ].

\par A ring is called a discrete valuation ring if it is a
principal ideal domain with only one maximal ideal. A ring $R$ is
a Dedekind ring if it is a Noetherian integral domain such that
the localization $R_{\fp}$ is discrete valuation ring for every
nonzero prime ideal $\fp$ of $R.$ Recall that a ring $R$ is a
Dedekind domain if and only if $R$ is Noetherian, integrally
closed domain and each nonzero prime ideal of $R$ is a maximal
ideal [\cite{J}, Theorem 3.16 p.13].
 \par Let $R$ be a Dedekind ring and let $I$ be a nonzero ideal of
 $R,$ we say that there exists some prime ideals $\fp_1,...,\fp_n$ (uniquely determined by
 $I$) and certain positive integers $k_1,...,k_n$ (uniquely determined by
 $I$) such that $I=\fp_1^{k_1}...\fp_n^{k_n}$ [\cite{J}, p.12].

 \par If $R$ is a ring and $E$ is an $R-$module, the idealization
 (also called trivial ring extension of $R$ by $E$) $R\propto E,$
 introduced by Nagata in 1956 is the set of pair $(r,e)$ with
 pairwise addition and multiplication given by
 $(r,e)(s,f)=(rs,rf+se).$ The trivial ring extension of $R$ by $E,\ R\propto E$
 has the following property that containing $R$ as sub-ring, where the module $E$ can be viewed as an ideal such that is square is zero.

\end{section}
\bigskip
\begin{section}{\ Main Results}
\bigskip
\begin{definition}
A ring $R$ is called a weak von Neumann regular ring (WVNR for
short) if for every finitely generated ideals $I$ and $J$ of $R$
satisfying that $I\subseteq J\varsubsetneq R,$ when $J$ is
generated by an idempotent element of $R,$ then so is $I.$
\end{definition}
In particular, any von Neumann regular ring is a weak von Neumann
regular ring. Now, we give a class of a weak von Neumann regular
ring.
\bigskip

\begin{proposition}\label{prop1}
Let $R$ be a ring. Assume that $R$ has exactly two idempotent
elements which are 0 and 1, then $R$ is a WVNR ring.
\end{proposition}
\begin{proof}
Let $(I,J)$ be a pair of proper finitely generated ideals such
that $I \subseteq J$ and $J$ is generated by an idempotent
element. Then $J=0,$ therefore $I=0.$ It follows that $R$ is a
WVNR ring.

\end{proof}
\bigskip
\begin{example} Let $R$ be an integral domain or a local ring. Then $R$ has
exactly two idempotent elements which are 0 and 1, so $R$ is a
WVNR ring.
\end{example}

\bigskip

 Now we give an example of a non-coherent WVNR ring.
\begin{example}
Let $X$ be a connected topological space and $R=\mathcal{C}(X,\R)$
the ring of numerical continuous functions defined in $X.$ Let $f$
be an idempotent element of $R,$ then for each $x\in X\; f(x)=0$
or $f(x)=1.$ Hence $$\forall x\in X\ f(x)=0 \; \mbox { or } \;
\forall x\in X\ f(x)=1$$ since $X$ is connected. Therefore $R$ has
exactly two idempotents and so $R$ is a WVNR ring.\\
Now we suppose that $X=[0,2].$ Let $f_0\in R$ such that $f_0(x)=0$
if $0 \leq x \leq 1$ and $f(x)\not=0$ if $1 < x \leq 2.$ Assume
that $(0:f_0)$ is a finitely generated ideal of $R,
(0:f_0)=(f_1,...,f_n).$ It is easy to see that $$(0:f_0)=\{\varphi
\in R\ ;\ \forall x\in [1,2]\ \varphi(x)=0\}.$$ Let $f\in R$
 defined by $f(x)=\sqrt{|f_1(x)|+...+|f_n(x)|}.$ Clearly we have
 $f\in(0:f_0),$ then there exists $(g_1,...,g_n)\in R^n$ such that
 $f=f_1g_1+...+f_ng_n.$ We claim that $$\forall r\in ]0,1[ \; \exists x\in[r,1]\ :\ f(x)\not=0
 .$$ Deny. There is some $r\in ]0,1[$ such that $f(x)=0$ for each $x\in
[r,2].$ Thus for any pair $(i,x)\in\{1,...,n\}\times [r,1[,$ we
have $f_i(x)=0.$
 Therefore $$(0:f_0)=\{\varphi
\in R\ ;\ \forall x\in [r,2]\ \varphi(x)=0\}.$$ Which is absurd.
 We conclude that for each nonnegative integer $p$
there exists $1-\frac{1}{p} \leq x_p \leq 1$ such that
$f(x_p)\not=0.$ On the other hand every $g_i$ is bounded mapping.
There is some $c>0$ such that $$\forall i\in \{1,...,n\} \;
\forall x\in [0,1] \; |g_i(x)|\leq c.$$ Thus $$f(x_p)\leq
c(|f_1(x_p)|+...+|f_n(x_p)|).$$ It follows that $1\leq cf(x_p)$ so
that in limit we get $\displaystyle 1\leq c\lim_{p \rightarrow
\infty }f(x_p).$ But $\displaystyle\lim_{p \rightarrow \infty
}f(x_p)=0,$ we have the desired contradiction. Consequently $R$ is
a non-coherent WVNR ring.
\end{example}

\bigskip

Now we give characterization of weak von Neumann regular rings.

\begin{theorem}\label{th1}
The following conditions on a ring $R$ are
equivalent:\begin{enumerate}
    \item $R$ is a WVNR ring.
    \item For each $a\in Re$ where $e$ is a nonunit idempotent
    element of $R,$ then $a\in Ra^2.$
    \item For each $a\in Re$ where $e$ is a nonunit idempotent
    element of $R,$ then $ Ra$ is a direct summand of $R.$
\end{enumerate}
\end{theorem}
\begin{proof}
(1) $\Rightarrow$ (3): Let $a\in R$ and $1\not= e$ an idempotent
element of $R$ such that $a\in Re.$ Since $e$ is nonunit we have
the containments $Ra\subseteq Re \subsetneq R.$ From the
definition of a WVNR ring, we can write $Ra=Rf$ for some
idempotent $f\in R.$
It follows that $Ra\oplus R(1-f)=R.$\\
\par
(3) $\Rightarrow$ (2): Let $a\in Re$ where $e$ is a nonunit
idempotent element of $R$ and let $I$ be an ideal of $R$ such that
$I\oplus Ra=R.$ We can write $1=u+v$ for some $u\in I$ and $v\in
Ra.$ Multiplying the above equality by $u$ (resp., $v$) we get
that $u^2=u$ (resp., $v^2=v$). Thus $I=Ru$ and $Ra=Rv,$ therefore
$a=au+av=av=a^2x$ for some $x\in R.$\\
\par
(2) $\Rightarrow$ (1): Let $J$ be a principal ideal generated by a
nonunit idempotent element $e$ of $R,$ and let $I$ be a finitely
generated ideal of $R$ contained in $J.$ It suffices to prove that
if $I=(a,b),$ then there exists an idempotent $f$ in $R$ such that
$I=Rf.$ Since $a\in J=Re$ then $a\in Ra^2,$ also $b\in Rb^2.$ Let
$u=ax$ and $v=by,$ where $a^2x=a$ and $b^2y=b.$ Hence $u$ and $v$
are idempotent elements of $R.$ The element $f=u+v-uv$ has the
required property.
\end{proof}
\bigskip
\begin{corollary}
Let $R$ be a ring. Then the following statements are
equivalent:\begin{enumerate}
    \item $R$ is a von Neumann regular ring.
    \item $R$ is a WVNR ring and for every nonunit element $a$ of
    $R$ there exists an idempotent $e\not=1$ of $R$ such that $a\in Re.$
\end{enumerate}
\end{corollary}

\begin{proof}
(1) $\Rightarrow$ (2): This implication is clear.\\
(2) $\Rightarrow$ (1): Let $a$ be a nonunit element of $R,$ assume
that there exists a nonunit idempotent element $e$ of $R$ such
that $a\in Re.$ By applying condition (2) of Theorem \ref {th1} we
get that $a\in Ra^2.$ Hence $a=a^2x$ for some $x\in R.$ Therefore
$R$ is a von Neumann regular ring.
\end{proof}

\bigskip

Now we give a necessary and sufficient condition for a direct
product of rings to be a WVNR ring.

\bigskip

\begin{theorem}\label{th2}
Let $(R_i)_{1\leq i\leq n}$ be a family of rings, with $n\geq2.$
Then the following statements are equivalent:\begin{enumerate}
    \item $\displaystyle\prod _{i=1}^n R_i$ is a von Neumann regular ring.
    \item $\displaystyle\prod _{i=1}^n R_i$ is a weak von Neumann regular ring.
    \item For each $i\in \{1,...,n\}\ R_i$ is a von Neumann regular
    ring.
\end{enumerate}
\end{theorem}
\begin{proof}
(1) $\Rightarrow$ (2): obvious.\\
(2) $\Rightarrow$ (3): Let $i\in \{1,...,n\},$ we need only show
that for each $r\in R_i$ there is some $s\in R$ such that
$r^2s=r.$ Let $e_i$ (resp., $a$) be the element of
$\displaystyle\prod _{j=1}^n R_j$ which has one (resp., $r$) in
its $i^{th}$ place and zeros elsewhere. It is easy to see that
$e_i$ is a nonunit idempotent element of $\displaystyle\prod
_{j=1}^n R_j$ and $a\in \displaystyle\left(\prod _{j=1}^n
R_j\right)e_i.$ Since $\displaystyle\prod _{j=1}^n R_j$ is a WVNR
ring, we have $\displaystyle a\in \left(\prod _{j=1}^n
R_j\right)a^2$ by using {Theorem}\ref{th1}. Hence there exists an
element $(x_1,...,x_n)\in \displaystyle\prod _{j=1}^n R_j$ such
that $a^2x=a,$ therefore $x_ir^2=r.$ It follows that $R_i$ is a
von Neumann
regular ring.\\
(3)$\Rightarrow$ (1): Let $\displaystyle a\in \prod _{j=1}^n R_j.$
For each $i\in\{1,...,n\},$ there exists an element $b_i\in R_i$
such that $a_i^2b_i=a_i.$ Then $a^2b=a,$ where $b=(b_1,...,b_n).$
Therefore $\displaystyle\prod _{j=1}^n R_j$ is a von Neumann
regular ring.
\end{proof}

\bigskip

Now we give an example of a non-WVNR Noetherian ring.

\begin{example}
Let $n$ be a positive integer such that $n\geq 2,$ then
$\mathbb{Z}^n$ is not a WVNR ring since $\mathbb{Z}$ is not a von
Neumann regular ring. Consequently, $\mathbb{Z}^n$ is non-WVNR
Noetherian ring.
\end{example}

\bigskip

For an ideal $I$ of a WVNR ring $R,\ R/I$ is not necessarily a
WVNR ring. For this, we claim that $\mathbb{Z}/12\mathbb{Z}$ is
not a WVNR ring. Indeed, 9 is an idempotent and $9.2=6$ but
$6^2=0.$ Thus $6\not=6^2x$ for each $x\in
\mathbb{Z}/12\mathbb{Z}.$ By applying condition (2) of
{Theorem}\ref{th1}, we get the result. \\

In the next theorem we give a sufficient condition for $R/I$ to be
a WVNR ring.

\bigskip

\begin{theorem}\label{th3}
Let $R$ be a ring and let $I$ be a primary ideal. Then $R/I$ is a
WVNR ring.
\end{theorem}

\begin{proof}
We denote $\overline{a}=a+I$ for every $a\in R.$ To prove
{Theorem}\ref{th3}, it is enough to show that $R/I$ has exactly
two idempotent elements which are $\overline{0}$ and
$\overline{1}.$ Let $a\in R$ such that $\overline{a}$ a nonzero
idempotent element of $R/I.$ We have $a^2-a\in I.$ Since $I$ is a
primary ideal of $R$ and $a\notin I,$ there exists a nonnegative
integer $n$ such that $(a-1)^n\in I.$ By the binomial theorem
(which is valid in any commutative ring),
$$(a-1)^n=\displaystyle\sum_{k=0}^{n}(-1)^{n-k}
\left(\begin{array}{c}
                                        n \\
                                         k \\
                                       \end{array}\right)a^k\in
                                       I.$$
 We put $a^2=a+x.$ By induction we claim that for
 each $k\geq 2,\ a^k=a+x\left(1+a+...+a^{k-2}\right).$
 Indeed, it is certainly true for $k=2.$  Suppose the statement is true
 for  $k,$ then we get the following equalities
 $$a^{k+1}=a^2+x\left(a+a^2+...+a^{k-1}\right)=a+x\left(1+a+...+a^{k-1}\right)$$
We conclude that for each $n\in \N^*,$ there is some $x_k\in I$
such that $a^k=a+x_k.$ We can also deduce that

 $$(-1)^n1+\displaystyle\sum_{k=1}^{n}(-1)^{n-k}
\left(\begin{array}{c}
                                        n \\
                                         k \\
                                       \end{array}\right)\left(a+x_k\right)\in
                                       I.$$
But$$(-1)^n1+\displaystyle\sum_{k=1}^{n}(-1)^{n-k}
\left(\begin{array}{c}
                                        n \\
                                         k \\
                                       \end{array}\right)a=(-1)^n(1-a),$$
hence $1-a\in I$ and so $\overline{a}=\overline{1}.$ By applying
{Proposition}\ref{prop1} we get that $R$ is a WVNR ring. This
completes the proof.
\end{proof}
\bigskip

\begin{remark}
The converse of {Theorem}\ref{th3} is not true in general. For
example $\mathbb{Z}/6\mathbb{Z}$ is a WVNR ring because
$$\forall x\in \mathbb{Z}/6\mathbb{Z} \  \; x^3-x=x(x-1)(x+1)=0.$$
But $6\mathbb{Z}$ is not a primary ideal of $\mathbb{Z}.$
\end{remark}

\bigskip

\begin{corollary}
Let $R$ be a Dedekind ring and let $I=\fp_1^{k_1}...\fp_n^{k_n}$
be a nonzero ideal of $R,$ where $\fp_1,...,\fp_n$ are the prime
ideals containing $I.$ Then $R/I$ is a WVNR ring if and only if
$n=1$ or $k_1=...=k_n=1.$
\end{corollary}

\begin{proof}
We shall need to use the following property:\\
if $\fp$ and $\fq$ are distinct maximal ideals of any ring $A$
then $\fp^k+\fq^l=A$ for every positive integers $k$ and $l.$
\par Assume that $n\geq 2.$ Thus $\fp_i^{k_i}+\fp_j^{k_j}=R$ if $i\not=j.$ By using
the Chinese remainder theorem we deduce that
$$R/I \simeq R/\fp_1^{k_1}\times...\times R/\fp_n^{k_n}.$$ We can now apply
{Theorem}\ref{th2} to obtain that $R/I$ is a WVNR ring if and only
if for each $i\in \{1,...,n\},R/\fp_i^{k_i}$ is a von Neumann
regular ring. On the other hand every power of a maximal ideal
$\fm$ of $R$ is a primary ideal. By applying {Theorem}\ref{th3} we
deduce that $R/\fp_i^{k_i}$ has exactly two idempotent elements
which are 0 and 1. But a WVNR ring $A$ is a von Neumann regular
ring if and only if for every nonunit element $a$ of $A$ there
exists a nonunit idempotent $e$ of $A$ such that $a\in Ae.$ It
follows that $R/I$ is a WVNR ring if and only for each
$i\in\{1,...,n\},\ R/\fp_i^{k_i}$ is a field (i.e
$\fp_i^{k_i}=\fp_i$).
\par Finally, if $n=1$ then $R/I$ is a WVNR ring by {Theorem}\ref{th3}. We conclude that $R/I$ is a WVNR ring if and only if $n=1$ or
$k_1=...=k_n=1.$
\end{proof}

\bigskip

\begin{example}
Let $n$ be a positive integer. Then $\mathbb{Z}/n\mathbb{Z}$ is a
WVNR ring if and only if $n$ is a power of a prime integer or
$v_p(n)\in \{0,1\}$ for every prime integer $p$ ($v_p(n)$ is the
$p-$valuation of $n$).
\end{example}

\bigskip

\begin{example}
Let $K$ be a field, and let $f$ be nonconstant polynomial in
$K[x]$. Thus $K[x]/(f)$ is a WVNR ring if and only if $f$ is a
power of a irreducible polynomial or $v_p(f)\in\{0,1\}$ for every
irreducible polynomial $p.$
\end{example}

\bigskip

Now, we give a characterization that a polynomial ring is a WVNR
ring.

\bigskip

\begin{theorem} \label{th4}
Let $R$ be a ring. Then the polynomial ring $R[x]$ is a WVNR ring
if and only if $R$ has exactly two idempotent elements.
\end{theorem}

\begin{proof}
Let $f=a_0+a_1x+...+a_nx^n$ an idempotent element of $R[x],$ then
$a_0^2=a_0.$ Let $k\in \{1,...,n\},$ by induction on k we prove
that $a_k=0.$ From the assumption, we can write $a_1=2a_0a_1.$
Multiplying this equality by $a_0$ we get that $a_1=0.$ We can
also write $\displaystyle\sum_{i=0}^ka_ia_{k-i}=a_k,$ then by
inductive hypothesis $a_k=2a_0a_k.$ Hence $a_k=0,$ therefore the
set of all idempotent elements of $R[x]$ is $\{a\in R\ ;\
a^2=a\}.$\\
\par Assume that $R$ has exactly two idempotent elements, from the
previous part of the proof we deduce that $R[x]$ is a WVNR ring.
Conversely, suppose that $R[x]$ is a WVNR ring and let $e$ be a
nonunit idempotent element of $R.$ We have $ex\in eR[x].$ By using
condition (2) of {Theorem}\ref{th1}, we get that $ex\in
(ex)^2R[x].$ There is some $f\in R[x]$ such that $ex=ex^2f(x).$
Thus $e=0,$ completing the proof of {Theorem}\ref{th4}.
\end{proof}

\bigskip

\begin{corollary}
Let $R$ be a ring. Then the polynomial ring $R[x_1,...,x_n]$ in
several indeterminates is a WVNR ring if and only if $R$ has
exactly two idempotent elements which are 0 and 1.
\end{corollary}

\begin{proof}
By induction on $n$ from {Theorem} \ref{th4}.
\end{proof}

\bigskip

\begin{example}
Let $R$ be a von Neumann regular ring which is not a field. Then
$R[x_1,...,x_n]$ is not a WVNR ring. For instance if $R$ is a
boolean ring such that $R\not= \{0,1\}$ then $R[x_1,...,x_n]$ is
not a WVNR ring.
\end{example}

\bigskip

\begin{remark}
Let $R$ be a ring and let $R[[x]]$ be the ring of formal power
series in $x$ with coefficients in $R.$ With a similar proof as in
{Theorem} \ref{th4}, we get that $R[[x]]$ is a WVNR ring if and
only if $R$ has exactly two idempotent elements.
\end{remark}

\bigskip
We end this paper by studying the transfer of a WVNR property to
trivial ring extensions.

\bigskip

\begin{theorem}\label{th5}
Let $A$ be a ring, $E$ an $A-$module and let $R=A\propto E$ be the
trivial ring extension of $A$ by $E.$ Then $R$ is a WVNR ring if
and only if the following statements are true:\begin{enumerate}
    \item $A$ is a WVNR ring.
    \item $aE=0$ for every idempotent element $1\not=a\in A.$
\end{enumerate}
\end{theorem}
\bigskip
Before proving { Theorem }\ref{th5}, we establish the following
lemma:

\begin{lemma}\label{lem2}
Let $A$ be a ring, $E$ an $A-$module and let $R=A\propto E$ be the
trivial ring extension of $A$ by $E.$ Then an element $(a,x)$ of
$R$ is idempotent if and only $a$ is idempotent and $x=0.$
\end{lemma}
\begin{proof}
Let $(a,x)$ be an idempotent element of $R,$ then
$\left(a^2,2ax\right)=(a,x).$ Hence $a^2=a$ and $2ax=x.$
Multiplying the second equality by $a$ we obtain that $x=0.$ We
deduce that $a$ is an idempotent element of $A$ and $x=0.$ The
converse is obvious.
\end{proof}
\bigskip
\begin{proof}$\textit{of \ theorem \ref{th5}}.$
Assume that $R$ is a WVNR ring. Let $a\in Ae$ for some nonunit
idempotent $e$ of $A,$ then $(a,0)\in R(e,0).$ The element $(e,0)$
is a nonunit idempotent of $R,$ by { Theorem }\ref{th1} we get
that $(a,0)\in R(a,0)^2.$ Hence there exists $(b,x)\in R$ such
that $(a,0)=(a,0)^2(b,x).$ Therefore $a\in Aa^2.$ We deduce that
$A$ is a WVNR ring. Now we consider a nonunit idempotent element
$a$ of $A$ and $x\in E.$ We have $(0,ax)=(a,0)(0,x),$ then
$(0,ax)\in R(a,0).$ Since $(a,0)$ is a nonunit idempotent element
of $R,$ then $\ (0,ax)\in R(0,ax)^2$ and so $ax=0.$ It follows
that $aE=0.$
\par Conversely, suppose that $A$ is a WVNR ring and $bE=0$ for
each nonunit idempotent element $b$ of $A.$ Let $(a,x)\in R(b,0)$
for some nonunit idempotent element $b$ of $A,$ there is some
$(c,y)\in R$ such $(a,x)=(b,0)(c,y).$ Hence $a\in Ab$ and $x=by$
therefore $x=0$ and $a\in Aa^2.$ It follows that $(a,x)\in
R(a,x)^2.$ This completes the proof of { Theorem }\ref{th5}.
\end{proof}
\bigskip
\begin{example}
Let $A$ be a ring and let $E$ be an $A-$module. Suppose that $A$
has exactly two idempotent elements. Then $A\propto E$ is a WVNR
ring. For instance let $G$ be a commutative group, then
$\mathbb{Z}\propto G$ is a WVNR ring. This is an other WVNR ring
which is neither local nor integral domain. Finally by { Theorem
}\ref{th4} the polynomial ring $\left(\mathbb{Z}\propto
G\right)[x]$ is also a WVNR ring.
\end{example}
\bigskip
\begin{corollary}
Let $A$ be a ring and $Q(A)$ its full ring of quotient. Then the
following statements are equivalent:\begin{enumerate}
    \item $A$ has exactly two idempotent elements.
    \item $A\propto A$ is a WVNR ring.
    \item $A\propto Q(A)$ is a WVNR ring.
\end{enumerate}
\end{corollary}
\begin{proof}
(1) $\Rightarrow$ (2): By { Lemma }\ref{lem2} $A\propto A$ has
exactly two idempotent elements which are 0 and 1.\\
(2) $\Rightarrow$ (1): Let $a$ a nonunit idempotent element of
$A.$ By { Theorem }\ref{th5} $aA=0,$ then $a=0.$\\
(1) $\Leftrightarrow $ (3): By the some way we get this
equivalence.
\end{proof}
\end{section}



\bigskip\bigskip


\end{document}